\documentclass[11pt]{amsart}
   \usepackage{amsmath,amssymb}
   
   \renewcommand{\bf}{\bfseries}

   \newcommand{\N}{\mathbb{N}}

   \newcommand{\R}{\mathbb{R}}

   \newtheorem{theorem}{Theorem}
   \newtheorem{lemma}{Lemma}
   
   \newtheorem{proposition}[theorem]{Proposition}
   \newtheorem{corollary}{Corollary}
   \renewcommand{\epsilon}{\varepsilon}

   \newcommand{\dis}{\displaystyle}
  
\begin{document}
  \title{On A. Zygmund differentiation conjecture}
   \author{I. Assani}
\thanks{Department of Mathematics,
UNC Chapel Hill, NC 27599, assani@email.unc.edu}
\thanks{AMS subject classification, 42B25, 37A05, 37A30}
\today
   \begin{abstract}
    Consider $v$ a Lipschitz unit vector field on $R^n$ and $K$ its Lipschitz constant. We show that the
    maps $S_s:S_s(X) = X + sv(X)$ are invertible for $0\leq |s|<1/K$ and define nonsingular point
    transformations.  We use these properties to prove first the differentiation in $L^p$ norm for
    $1\le p<\infty.$ Then we show the existence of a universal set of values $s\in [-1/2K,1/2K]$ of
    measure $1/K$ for which the Lipschitz unit vector fields $v\circ S_s^{-1}$ satisfy Zygmund's conjecture
    for all functions in $L^p(\R^n)$ and for each p,  $1\leq p< \infty.$
  \end{abstract}
  \maketitle
  \section{Introduction}
Lebesgue differentiation theorem states that given a function $f\in
L^1(\R)$ the averages $\dis \frac{1}{2t}\int_{-t}^t f(x+u) du$
converge a.e. to $f(x)$ when $t$ tends to zero. The differentiation
for functions $F$ defined on $\R^2$ is more subtle. Actually it is a
longstanding problem to find analogue of Lebesgue differentiation
theorem for averages of the form $$ M_t(F)(x,y)=
\frac{1}{2t}\int_{-t}^t F[(x,y)+ \beta v(x,y)]d\beta$$ for a
measurable function $v$. One would expect these averages to converge
a.e. to $F(x,y).$ In other words one looks at the differentiation
along the vector field $v$ (or the direction $v$).(see for instance
\cite{Stein}, \cite{DeGuzman}). One can see that because of the
geometry of $R^2$ multiple directions are possible. In fact the
example of the Nikodym set \cite{DeGuzman} shows that condition on
$v$ must be imposed if one expects the differentiation to hold. J.
Bourgain \cite{Bourgain} established the differentiation of the
averages $M_t(F)$ for function $F\in L^2$ and $v$ a real analytic
vector field. N.H. Katz \cite{Katz} has some partial result for
Lipschitz vector fields. A longstanding conjecture attributed to A.
Zygmund (see the paper by M. Lacey and X. Li,  \cite{LaceyLi}) is
the following. \vskip1ex
 \noindent{\bf {Zygmund's conjecture}}
\vskip1ex {\em Let $v$ be a Lipschitz unit vector field and let
$F\in L^2(\R^2)$. Do the averages
  $$ M_t(F)(x,y)=\frac{1}{2t}\int_{-t}^t F[(x,y)+ \beta v(x,y)]d\beta$$ converge
a.e. to $F(x,y)$?}

\vskip1ex
First we will observe that for $s$ small enough (if $K$ is the
Lipschitz constant of $v$ we will require $|s|< 1/2K$), the maps
$\dis S_s: S_s(x)= x + sv(x)$ are invertible. This observation will
allow us to derive the norm convergence of the averages
$$M_t(F)= \frac{1}{2t}\int_{-t}^t F(x+ \beta v(x))d\beta$$ to the function F
in all $L^p$ spaces, $1\leq p <\infty.$ This norm convergence result
was apparently an open problem (see \cite{Bourgain}.)

Then we will show that Zygmund's conjecture holds in all $L^p$
spaces $1\leq p <\infty$ for the unit vector fields $v\circ
S_s^{-1}$ when $s\in \mathcal{T},$ a universal subset
 of  $[-1/2K, 1/2K]$ with measure $1/K.$
 The method we use extends to $R^n.$
 \vskip1ex
 \noindent{\bf Acknowledgments} We thank C. Thiele and C. Demeter for
bringing this problem to our attention. Thanks also to C. Demeter
for his comments on a preliminary version of the paper.
\vskip1ex

\section{Differentiation in $\R^2$}
 The main steps are as follows. First we show that
 for $s$ small enough the maps $S_s:S_s(x,y) = (x, y) + sv(x,y)$
 are invertible and nonsingular in the sense that $\mu(A)=0$ if
 and only $\mu(S_s(A))= 0.$ A more precise statement is given in
 Lemma 1 where we prove that the operators induced by these maps
 are uniformly bounded on $L^p(\R^2)$ for $1\leq p\leq \infty.$
  From this we derive the norm convergence of the averages $M_t(F)$
 to $F$. Two consequences are derived from Lemma 1. First we obtain a "weak" version of our main result,
 Proposition 2, where we show that given a function $F\in L^1(\R^2)$ the
differentiation occurs along the vector fields $v\circ S_s^{-1}$
 as long as $s$ belongs to a set of measure $1/K$ depending a priori on $F.$ Then we use Hardy Littlewood maximal
 inequality on $L^1(\R)$ to derive a first maximal inequality for the
differentiation problem (Theorem 3).  Our main result is proved by
showing that the set where the differentiation occurs can in fact be
taken independently of any $F\in L^1(\R^2).$ Finally we establish a
"local" maximal inequality for the maximal operator associated with
these averages.
\subsection{Convergence in $L^p$ norm}
 \begin{lemma}
Assume that $v$ is a unit vector field (i.e $\|v(x,y)\|_2= 1$ for
all $(x,y) \in R^2$ and a Lipschitz map with constant $K$. Then for
 each $t; |t|\leq T<\frac{1}{K}$ the map $S_t$ from $R^2$ to $R^2$
such that $S_t(x, y) = (x,y) + tv(x,y)$ is one to one and onto.
Furthermore if we denote by $\mu$ Lebesgue measure on $R^2$ for all
measurable sets $A\subset R^2,$ for all $|s|\leq T,$ we have
$$\frac{1}{2\pi(1+ |s|K)^2}\mu(S_s(A)) \leq \mu(A)\leq 2\pi
\big(\frac{1}{1-|s|K}\big)^2\mu(S_s(A)).$$
\end{lemma}
\begin{proof}

\vskip1ex
 First if $S_t(x_1,y_1)= S_t(x_2, y_2)$ then we have
$$\|(x_1,y_1)- (x_2, y_2)\|= \|t(v(x_1, y_1)- v(x_2, y_2)\|\leq KT \|(x_1,y_1)- (x_2, y_2)\|.$$
As $KT<1$ this shows that $S_t$ is one to one. \vskip1ex
 The equation $Z=(z_1, z_2)= (x,y) + tv(x,y)= X + tv(X)$ has a solution in
 $X=(x,y)$ that can be found by applying the fixed point theorem to
 the function $R_Z$:$R_{Z}(X)= Z + X - S_t(X).$
 \vskip1ex
 To establish the second part of the lemma we can observe that it
 is enough to prove it for cubes $A.$
 For any two points $Z_1= X_1 + sv(X_1)$ and $Z_2= X_2 + sv(X_2)$
 we have $\|Z_1- Z_2\|\leq (1+ |s|K)\|X_1- X_2\|,$ and
 $\|X_1-X_2\|\leq \frac{1}{1-|s|K} \|S_s(X_1)-S_s(X_2)\|.$
 Also for each measurable compact set  $B\subset R^2$ we have $\mu (B)
 \leq \pi diam(B)^2,$ where $diam(B)$ is the diameter of the bounded set $B.$
 As $\|S_s(X)- S_s(Y)\| \leq (1+ |s|K)\|X-Y\|$ we have
 $diam(S_s(A)) \leq (1+|s|K)diam(A).$ Therefore if we denote by $r$ the side length of the cube $A$ we have
 $$\mu (S_s(A))\leq \pi diam(S_s(A))^2\leq \pi (1+|s|K)^2 diam(A)^2=
 \pi (1+|s|K)^2 2r^2\leq 2\pi (1+|s|K)^2\mu(A).$$
 By approximation we conclude that for any measurable set $A$ we
 have the same inequality.
 From the inequality
$\|X_1-X_2\|\leq \frac{1}{1-|s|K} \|S_s(X_1)-S_s(X_2)\|,$
 we can conclude that $$\|S_s^{-1}(Y_1)-
 S_S^{-1}(Y_2)\|
\leq \frac{1}{1-|s|K}\|Y_1- Y_2\|$$ for all $Y_1, Y_2\in R^2.$ The
same path will lead us then to the inequality
$$\mu(S_s^{-1}(B)) \leq
2\pi \big(\frac{1}{1-|s|K}\big)^2 \mu (B)$$ for all measurable set
$B\subset R^2.$ From this we can derive the second inequality in the
lemma.
\end{proof}

\vskip1ex Using the notations of Lemma 1 we can obtain the
convergence in $L^p$ norm.
\begin{proposition}
For $0<|t|\leq T$ and for $1\leq p \leq \infty$ the operators $M_t$
defined pointwise by
 $$M_t(F)(x,y) = \frac{1}{2t}\int_{-t}^t F[(x,y) + sv(x,y)]ds
$$ map $L^p$ into $L^p.$ Furthermore for each $1\leq p<\infty$ we
have
$$\lim_{t\rightarrow 0}\|M_t(F)- F\|_p = 0.$$
\end{proposition}
\begin{proof}
 It follows immediately from Lemma 1. Indeed the case $p= \infty$ is obvious.
 For the other values of $p$, consider
 a nonnegative simple $L^p$ integrable function $\dis F= \sum_{n=1}^N \alpha_n\mathbf{1}_{A_n}$ with disjoint measurable
  sets $A_n$. We have
 \[
 \begin{aligned}
 &\|M_t(F)\|_p^p =
 \int_{\R}\bigg|\frac{1}{2t}\int_{-t}^{t}\sum_{n=1}^N
 \alpha_n\mathbf{1}_{A_n}(S_s(x,y))ds\bigg|^p d\mu\\
 &\leq \int_{\R}\bigg(\frac{1}{2t}\int_{-t}^{t}\sum_{n=1}^N
 \alpha_n^p\mathbf{1}_{A_n}(S_s(x,y)\bigg)ds d\mu =
 \frac{1}{2t}\int_{-t}^{t}\sum_{n=1}^N
 \alpha_n^p\mu(S_s^{-1}(A_n))ds\\
 &\leq \frac{1}{2t}\bigg(\int_{-t}^t 2\pi\big(\frac{1}{1-|s|K}\big)^2 ds\bigg)\sum_{n=1}^N
 \alpha_n^p\mu(A_n) = \frac{2\pi}{1-tK}\|F\|_p^p
 \end{aligned}
 \]
  The boundedness of the operators $M_t$ follows by approximation.
  \vskip1ex
  The second part of the proposition is a consequence of the simple
  fact that for the dense set of continuous functions with compact support we have
  the pointwise and norm convergence of the operators $M_t.$
\end{proof}
\subsection{A "weak" version of Zygmund's conjecture}
\vskip1ex
 The next proposition is a "weak" version of Zygmund's
 conjecture in the sense that for each function $F\in L^1(\mu)$
 there exists a set of $s$ of measure $T$ such that
 $$\lim_{t\rightarrow 0}\frac{1}{2t}\int_{-t}^t F[(x,y) +
 \beta v(S_s^{-1}x)d\beta= F(x,y)$$ for almost every $(x,y)\in
 R^2.$ In other words the set of $s$ and Lipschitz vector fields $v\circ S_s^{-1}$ for which the differentiation
 occurs may depend on $F.$
 The next proposition gives us also a path on how to approach Zygmund's conjecture, more precisely by considering the averages along the values of the function $F$ at
 $(x,y) + \beta v(S_s^{-1}(x,y))$ and by exploiting the invertibility of
 the maps $S_s.$
  \begin{proposition}
  Let $v$ be a Lipschitz function from $\R^2$ to $R^2$ with Lipschitz
  constant $K$ such that $\|v(x,y)\|_2=1$ for almost all $(x,y)\in \R^2.$
   Then for all function $F\in L^1(R^2)$ for almost every $s\in [-T/2, T/2]$,
   for almost every $(x,y) \in R^2$ we have
   $$\lim_{t\rightarrow 0} \frac{1}{2t}\int_{-t}^t F[(x,y)+ (s+ \beta)v(x,y) ]d\beta =
   F[(x,y)+ sv(x, y)]$$
   and
   $$\lim_{t\rightarrow 0} \frac{1}{2t}\int_{-t}^t F[(x,y)+ \beta v(S_{s}^{-1}(x,y)]d\beta = F(x,y).$$
  \end{proposition}


\begin{proof}
 For $t, s$ and $\beta$ small enough we
consider the averages
$$\frac{1}{2t}\int_{-t}^{t} F[(x,y)+ (s+\beta)v(x,y)]d\beta$$
Because of the assumptions made on $v$ by Lemma 1
   for each $s; |s|\leq T<\frac{1}{K}$ for almost all $(x, y)\in R^2$
   $$G_{x,y}(s)=  F[(x,y)+ s v(x,y)] $$ is well defined and $G_{x,y}\in L^1([-T,T]).$
   By Lebesgue differentiation theorem, for almost every $s\in [-T/2, T/2]$ we have
   $$\lim_t \frac{1}{2t}\int_{-t}^t F[(x,y)+ (s+ \beta) v(x,y)]d\beta = F[(x,y)+ sv(x,y)].$$
   Let us consider the complement $E$ in $\R^2\times [-T/2, T/2]$ of the set
   $$ \{ (x, y, s): \lim_t \frac{1}{2t}\int_{-t}^t F[(x,y)+ (s+ \beta) v(x,y)]d\beta = F[(x,y)+
   sv(x,y)]\}$$
   By Fubini this set has measure zero. Again by Fubini for almost
   all $s$ the set
   $E_s = \{(x, y): (x,y, s)\in E\}$ also has measure zero.
   By lemma 1 the corresponding sets $S_s(E_s)$ will also have
   measure zero. This proves the second part of the proposition.
   \vskip1ex
   \end{proof}
   \vskip1ex
    As indicated above the maximal inequality allowing to derive
    the conclusions of proposition 3 is given by the following
    result.

\begin{theorem}
 Let $K$ be the Lipschitz constant for the unit vector field $v.$ Then for each $T$, $0<T< 1/K,$ for all $\lambda>0$
\[
\begin{aligned}
&\frac{1}{T}\int_{-T/2}^{T/2}\mu\bigg\{(x,y)\in R^2: \sup_{0<t\leq
T/2}\frac{1}{2t}\int_{-t}^t |F[(x,y) + \beta
v(S_s^{-1}(x,y))]|d\beta >\lambda\bigg\}dm(s)\\
&\leq \frac{4\pi^2 (1+ TK)^2}{(1-TK)^2}\frac{1}{\lambda}\int_{\R^2}
|F(x,y)|d\mu.
\end{aligned}
\]
 where $m$ denotes Lebesgue measure on $[-T/2, T/2].$
\end{theorem}
\begin{proof}
\vskip1ex
 For a.e. $(x, y)$ the function $G(x,y): G_{x,y}(s)=
\mathbf{1}_{[-T , T]}(s). F[(x,y)+ s v(x,y)] $ belongs to $L^1.$ By
Hardy-Littlewood maximal inequality applied to this function
 we have;
\[
\begin{aligned}
& m\bigg\{s\in [-T/2, T/2]: \sup_{0<t\leq
T/2}\frac{1}{2t}\int_{-t}^t |F[(x,y)+ (s + \beta )v(x,y)]| d\beta > \lambda \bigg\} \\
&\leq \frac{1}{\lambda} \int_{-T}^{T} |F[(x,y)+ \beta
v(x,y)]|d\beta.
\end{aligned}
\]

 We can integrate both sides of this inequality with respect to
Lebesgue measure $\mu$ on $\R^2$ and apply Fubini theorem.

We obtain by using Lemma 1,
\[
\begin{aligned}
&\mu\times{m} \bigg\{(x,y,s)\in R^2\times [-T/2, T/2] :
\sup_{0<t\leq T/2}\frac{1}{2t}\int_{-t}^t |F[(x,y)+ (s+ \beta) v(x,y)]|d\beta > \lambda \bigg\} \\
&\leq \frac{1}{\lambda} \int_{-T}^{T}\int_{\R^2}
|F[(x,y)+ \beta v(x,y)]|d\mu d\beta \\
&\leq \frac{2\pi}{(1-TK)^2\lambda}\int_{-T}^{T}\int_{\R^2}|F(x,y)|d\mu d\beta \\
&= \frac{2\pi T}{(1-TK)^2\lambda}\int_{\R^2} |F(x,y)|d\mu =
\frac{CT}{\lambda} \int_{\R^2} |F(x,y)|d\mu
\end{aligned}
\]
Dividing all expressions above by $T$ and rewriting
$$\mu\times{m}\bigg\{(x,y,s)\in R^2\times [-T/2, T/2] :  \sup_{0<t\leq
T/2}\frac{1}{2t}\int_{-t}^t |F[(x,y)+ (s+ \beta) v(x,y)]|d\beta
 > \lambda \bigg\}
 $$
 as
 $$\int_{-T/2}^{T/2} \mu\bigg\{(x,y)\in R^2:\sup_{0<t\leq
T/2}\frac{1}{2t}\int_{-t}^t |F[(x,y)+ (s+ \beta) v(x,y)]| d\beta
 > \lambda \bigg\}dm(s),$$
 we derive the following inequality:
 \[
 \begin{aligned}
&\frac{1}{T}\int_{-T/2}^{T/2} \mu\bigg\{(x,y)\in R^2:\sup_{0<t\leq
T/2}\frac{1}{2t}\int_{-t}^t |F[(x,y)+ (s+ \beta) v(x,y)]|d\beta
 >\lambda \bigg\}dm(s) \\
 &\leq \frac{C}{\lambda} \int_{\R^2} |F(x,y)|d\mu.
 \end{aligned}
 \]

 Using Lemma 1 we can observe that
 \[
 \begin{aligned}
 &\mu\bigg\{(x,y)\in R^2:\sup_{0<t\leq
T/2}\frac{1}{2t}\int_{-t}^t |F[(x,y)+ (s+ \beta) v(x,y)]|d\beta >\lambda\bigg\}\\
 &\geq \frac{1}{2\pi(1+ TK)^2}\mu\bigg\{(x,y)\in R^2: \sup_{0<t\leq
T/2}\frac{1}{2t}\int_{-t}^t |F[(x,y)+ \beta v(S_s^{-1}(x,y))]|d\beta
>\lambda\bigg\}.
 \end{aligned}
 \]
Therefore, for all $\lambda>0,$  we have the inequality
\[
\begin{aligned}
&\frac{1}{T}\int_{-T/2}^{T/2}\mu\bigg\{(x,y)\in R^2: \sup_{0<t\leq
T/2}\frac{1}{2t}\int_{-t}^t |F[(x,y) + \beta
v(S_s^{-1}(x,y))]|d\beta >\lambda\bigg\}dm(s)\\
&\leq \frac{4\pi^2 (1+ TK)^2}{(1-TK)^2}\frac{1}{\lambda}\int_{\R^2}
|F(x,y)|d\mu.
\end{aligned}
\]

\end{proof}
\subsection{A universal set of unit Lipschtiz vector fields
satisfying Zygmund's conjecture in all $L^p$ spaces.}

\vskip1ex

  As indicated in the introduction we want to strengthen Proposition 2 by showing that a
universal set of vector fields $v\circ S_s^{-1}$ satisfy Zygmund's
conjecture. More precisely we want to prove the following result.

\begin{theorem}
Let $v$ be a unit Lipschitz vector field with Lipschitz constant $K=
1/T$. Then there exists a set $\mathcal{T}\subset [-T/2, T/2]$ of
measure $T$ such that for each $s\in \mathcal{T}$ the unit Lipschitz
vector field $v\circ S_s^{-1}$ satisfies Zygmund's conjecture in all
$L^p$ spaces for $1\leq p<\infty$. More precisely for all $F\in
L^p(\R^2)$ the averages
$$\frac{1}{2t}\int_{-t}^t F[(x,y) + \beta v(S_s^{-1}(x,y))d\beta$$
converge a.e to $F(x,y).$
\end{theorem}

To prove this theorem we introduce some notation. We denote by $\dis
\mathcal{D}_N = \{(x,y)\in \R^2: \|(x,y)\|\leq N\}$ and by
$\mathcal{E}$ a countable set of continuous functions with compact
support dense in the unit closed ball of $L^1(\R^2).$ We will use
the notation $\dis M_t^s(F)(x,y)$ for the averages
$$\frac{1}{2t}\int_{-t}^t |F[(x,y) + \beta v(S_s^{-1}(x,y))]|d\beta$$

 \vskip1ex
Theorem 4 is a consequence of the following result.
\begin{theorem}
Under the assumptions of Theorem 4, we have for each p, $ 1\leq
p<\infty$ and for a.e. $s\in [-T/2, T/2]$
$$ \lim_{n\rightarrow \infty}\sup_{\|F\|_p\leq 1}\mu\{(x,y)\in
\mathcal{D}_N: \sup_{0<t<T/2}M_t^s(F)(x,y)> n\} = 0.$$
\end{theorem}
Our proof of these theorems will require several lemmas. We only
give the proof for the case $p=1.$ The case $p>1$ can be obtained
similarly without difficulty as the differentiation is a local
property. Given any function $F\in L^1(\mu)$ there exists a
subsequence $\dis G_j= F_{N_j}$ such that $\lim_j G_j(x,y) = F(x,y)$
except on a set of measure zero $\mathcal{N}$. The next lemma is a
consequence of Lemma 1. It shows that for almost every $(x,y)$ we
can keep this convergence along the line segments $[(x,y)- \beta
w(x,y), (x,y) + \beta w(x,y)]]$ where $\beta$ is in absolute value
smaller that the reciprocal of the Lipschitz constant of the unit
vector field $w$.

\begin{lemma}
Let $F\in L^1(\mu)$ and $G_j$ a sequence of continuous function with
compact support converging a.e. to $F$ . Let $w$ be a unit vector
field with Lipschitz constant {\bf{K}} . Then for almost all $x\in
\R^2$ for all $t\in [-T/2, T/2]$ we have
\[
\begin{aligned}
&\frac{1}{2t}\int_{-t}^t |F[(x,y)+ \beta w(x,y)]|d\beta=
\frac{1}{2t}\int_{-t}^t \liminf_j|G_j|[(x,y)+ \beta w(x,y)]d\beta \\
& = \sup_j\frac{1}{2t}\int_{-t}^t \inf_{k\geq j} |G_k|[(x,y)+ \beta
w(x,y)]d\beta
\end{aligned}
\]
\end{lemma}
\begin{proof}
Let us consider the null set $\mathcal{N}$ off which the sequence
$G_j$ converges to $F$. We can assume that this set is measurable.
Hence by Fubini we have
\[
\begin{aligned}
&\int_{\R^2}\int_{-T}^T \mathbf{1}_{\mathcal{N}}[(x, y) + \beta
w(x,y)]d\beta d\mu = \int_{-T}^T
\int_{\R^2}\mathbf{1}_{\mathcal{N}}[(x, y) + \beta w(x,y)]d\beta
d\mu \\
&\leq \big(\frac{1}{1-TK}\big)^2 2T\mu (\mathcal{N})= 0.
\end{aligned}
\]
Therefore there exists a set $A$ of zero measure such that for
$(x,y)\in A^c$ we have \vskip1ex $\dis \int_{-T}^T
\mathbf{1}_{\mathcal{N}}[(x, y) + \beta w(x,y)]d\beta= 0$. Hence for

all $t\in [-T/2, T/2]$ we also have $\dis \int_{-t}^t
\mathbf{1}_{\mathcal{N}}[(x, y) + \beta w(x,y)]d\beta= 0.$
 Writing the function $|F|$ as $\mathbf{1}_{\mathcal{N}}|F| +
 \mathbf{1}_{\mathcal{N}^c}|F|$ we have then for $(x, y)\in A^c$
 \[
 \begin{aligned}
 & \int_{-t}^t |F[(x, y) + \beta w(x,y)]|d\beta = \int_{-t}^t
 \mathbf{1}_{\mathcal{N}^c}[(x,y) + \beta w(x, y)]|F[(x,y) + \beta
 w(x,y)]|d\beta\\
 & = \int_{-t}^t\mathbf{1}_{\mathcal{N}^c}[(x,y) + \beta w(x, y)]\liminf_j|G_j|[(x,y)+ \beta w(x,y)]d\beta \\
& = \lim_j\frac{1}{2t}\int_{-t}^t \mathbf{1}_{\mathcal{N}^c}[(x,y) +
\beta w(x, y)]\inf_{k\geq j} |G_k|[(x,y)+ \beta w(x,y)]d\beta \\
& \text{by the monotone convergence theorem} \\
& = \lim_j\frac{1}{2t}\int_{-t}^t \inf_{k\geq j} |G_k|[(x,y)+ \beta
w(x,y)]d\beta
\end{aligned}
\]

\end{proof}
Next we want to check that the preceding lemma applies to all
Lipschitz unit vector fields $v\circ S_s^{-1}.$

\begin{lemma}
Let $v$ be a Lipschitz unit vector field with Lipschitz constant
$K$. Consider $0<T<1/K.$ Then for all $s$ such that $0<|s|< T/2$ the
unit vector fields $v\circ S_s^{-1}$ are Lipschitz vector fields
with Lipschitz constant $2K.$
\end{lemma}
\begin{proof}
As $v$ is a Lipschitz vector field with Lipschitz constant $K$ we
 have for all $X$, $Y$ in $\R^2$
$$\|v(S_s^{-1})(X) - v(S_s^{-1})(Y)|\leq K\|S_s^{-1}(X) -
S_s^{-1}(Y)\|.$$

We denote $Z_1= S_s^{-1}(X)$ and $Z_2 = S_s^{-1}(Y).$ Then
 $\dis X = Z_1 + sv(Z_1)$ and $\dis Y = Z_2 + sv(Z_2).$
 Therefore
 $\dis \|Z_1- Z_2\|\leq \|X - Y\| + |s|K\|Z_2- Z_1\|$ and we obtain
 \vskip1ex

 $\dis \|Z_1 - Z_2\|\leq \frac{1}{(1-|s|K)}\|X-Y\|.$
We conclude then that
 $$\|v(S_s^{-1})(X) - v(S_s^{-1})(Y)|\leq K\frac{1}{(1-|s|K)}\|X-Y\|.$$
 Noticing that for $0<|s|<T/2$ we  have $\dis \frac{1}{(1-|s|K)}\leq
 2$ and this concludes the proof of this lemma.
\end{proof}

\vskip1ex Thus we can apply Lemma 2 with the constant $\bold{K}=
2K.$

\begin{lemma}
For each $\lambda>0$ and each $s\in [-T/2, T/2]$ we have

\begin{equation}
\begin{aligned}
 & \sup_{\|F\|_1\leq 1}\mu\bigg\{(x,y)\in \mathcal{D}_N:
 \sup_{0<t<T/2}\frac{1}{2t}\int_{-t}^t |F[(x,y) + \beta
v(S_s^{-1}(x,y)]|d\beta >\lambda \bigg\} \\
&=\sup_{\Phi_j\in \mathcal{E}}\mu\bigg\{(x,y)\in \mathcal{D}_N:
 \sup_{0<t<T/2}\frac{1}{2t}\int_{-t}^t |\Phi_j[(x,y) + \beta
v(S_s^{-1}(x,y)]|d\beta >\lambda \bigg\}
\end{aligned}
\end{equation}
\end{lemma}
\begin{proof}
Let us fix $\epsilon>0$. We can find a function $F\in L^1$ with
$\|F\|_1\leq 1$ such that
\begin{equation}
\begin{aligned}
& \sup_{\|G\|_1\leq 1}\mu\bigg\{(x,y)\in \mathcal{D}_N:
 \sup_{0<t<T/2}\frac{1}{2t}\int_{-t}^t |G[(x,y) + \beta
v(S_s^{-1}(x,y)]|d\beta >\lambda \bigg\} \\
&\leq \mu\bigg\{(x,y)\in \mathcal{D}_N:
\sup_{0<t<T/2}\frac{1}{2t}\int_{-t}^t |F[(x,y) + \beta
v(S_s^{-1}(x,y)]|d\beta >\lambda \bigg\} + \epsilon
\end{aligned}
\end{equation}

 For the function $F$ we can find a subsequence $G_j= F_{n_j}$ of continuous
 functions in $\mathcal{E}$ which converges a.e. to $F.$
 Applying Lemma 2 with $w= v\circ S_s^{-1}$ off a null set
 $\mathcal{N}_s$ we have
 \[
 \begin{aligned}
&\sup_{0<t<T/2}\frac{1}{2t}\int_{-t}^t |F[(x,y)+ \beta
w(x,y)]|d\beta=
\sup_{0<t<T/2}\frac{1}{2t}\int_{-t}^t \liminf_j|G_j|[(x,y)+ \beta w(x,y)]d\beta \\
& = \sup_{0<t<T/2}\sup_j\frac{1}{2t}\int_{-t}^t \inf_{k\geq j}
|G_k|[(x,y)+ \beta w(x,y)]d\beta\\
& = \sup_{j}\sup_{0<t<T/2}\frac{1}{2t}\int_{-t}^t \inf_{k\geq j}
|G_k|[(x,y)+ \beta w(x,y)]d\beta \\
& = \lim_j \sup_{0<t<T/2}\frac{1}{2t}\int_{-t}^t \inf_{k\geq j}
|G_k|[(x,y)+ \beta w(x,y)]d\beta \\
& \text{(Noticing that the sup is the limit because we have an
increasing sequence)}
\end{aligned}
\]

 Hence we have
\[
\begin{aligned}
&\mu\bigg\{(x,y)\in \mathcal{D}_N:
 \sup_{0<t<T/2}\frac{1}{2t}\int_{-t}^t |F[(x,y) + \beta
v(S_s^{-1}(x,y)]|d\beta >\lambda \bigg\} \\
&= \mu\bigg\{(x,y)\in \mathcal{D}_N: \lim_j
\sup_{0<t<T/2}\frac{1}{2t}\int_{-t}^t \inf_{k\geq j} |G_k|[(x,y)+
\beta v(S_s^{-1}(x,y)]d\beta >\lambda\bigg\} \\
&= \lim_j\mu\bigg\{(x,y)\in \mathcal{D}_N:
\sup_{0<t<T/2}\frac{1}{2t}\int_{-t}^t \inf_{k\geq j} |G_k|[(x,y)+
\beta v(S_s^{-1}(x,y)]d\beta>\lambda\bigg\}\\
&\leq \limsup_j\mu\bigg\{(x,y)\in \mathcal{D}_N:
\sup_{0<t<T/2}\frac{1}{2t}\int_{-t}^t \inf_{k\geq j} |G_k|[(x,y)+
\beta v(S_s^{-1}(x,y)]d\beta>\lambda\bigg\} \\
&\leq\sup_{\Phi\in\mathcal{E}}\mu\bigg\{(x,y)\in
\mathcal{D}_N:\sup_{0<t<T/2}\frac{1}{2t}\int_{-t}^t |\Phi|[(x,y)+
\beta v(S_s^{-1}(x,y)]d\beta >\lambda\bigg\}
\end{aligned}
\]
 This last inequality combined with (2) proves Lemma 4.
\end{proof}

\begin{lemma}
For $F$ continuous with compact support and each $\lambda>0$ the map
$$s\in [-T/2, T/2]:\rightarrow \mu\bigg\{(x,y)\in \mathcal{D}_N;
\sup_{0<t<T/2} M_t^s(F)(x,y)>\lambda \bigg\}$$ is continuous.
\end{lemma}
\begin{proof}
Again we denote by $X$ the vector  $(x,y)\in \R^2$. For all
$|\beta|\leq T/2$  and for all $|s_1|, |s_2|\leq T/2$ we have
$$\|(X + \beta v(S_{s_1}^{-1}X)) - (X + \beta
v(S_{s_2}^{-1}X))\|\leq \frac{T}{2}K\|S_{s_1}(X) - S_{s_2}(X)\|.$$
For $Z_1 = S_{s_1}(X)$ and $Z_2= S_{s_2}(X)$ we have
$$Z_1 + s_1 v(Z_1)= X = Z_2 + s_2v(Z_2).$$
Therefore we have
$$Z_1 - Z_2 = s_2 v(Z_2) - s_1v(Z_1) = (s_2-s_1)v(Z_2) +
s_1(v(Z_2)-v(Z_1)).$$ As a consequence we obtain
$$\|Z_1- Z_2\| \leq |s_2- s_1| + \frac{T}{2}K \|Z_1- Z_2\|$$
and this gives us the uniform estimate
$$\|(X + \beta v(S_{s_1}^{-1}X)) - (X + \beta
v(S_{s_2}^{-1}X))\|\leq \frac{KT}{2}{(1- \frac{TK}{2})}|s_1- s_2|=
C|s_1- s_2|.$$ Now we can conclude by using the uniform continuity
of the function $F.$ For $|s_1-s_2|<\frac{\delta(\epsilon)}{C}$ then
for all $X\in R^2$ we have
$$\bigg|\sup_{0<t<T/2}M_t^{s_1}F(X) -
\sup_{0<t<T/2}M_t^{s_2}F(X)\bigg|<\epsilon.$$

\end{proof}

 The following lemma is well known and can be found in \cite{Folland}. We just
 state it to make the paper hopefully easier to read.
 \begin{lemma}
Let $\mathcal{C}$ be any collection of open intervals $B$ in $\R$
and let $U$ be the union of all these open intervals. If $c< m(U),$
then there exist disjoint $B_1, ...,B_k\in \mathcal{C}$ such that
$\dis \sum_{j=1}^k m(B_j)>\frac{1}{3} c.$
 \end{lemma}

Now we can proceed with the proof of Theorem 5. For simplicity we
will denote by $M_*^s(F)(X)$ the maximal function
$$\sup_{0<t<T/2}M_t^s(F)(x,y)$$

\vskip1ex \noindent{\bf{Proof of Theorem 5}}

\vskip1ex

 We will argue by contradiction.  Because of Lemma 4
 the functions
$$H_n: s\in [-T/2, T/2]\rightarrow H_n(s)=\sup_{\|F\|_1\leq 1}\mu\bigg\{(x,y)\in
\mathcal{D}_N: \sup_{0<t<T/2}M_t^s(F)(x,y)> n\bigg\}$$ being equal
for each s to
$$\sup_{F_i\in \mathcal{E}}\mu\bigg\{(x,y)\in
\mathcal{D}_N: \sup_{0<t<T/2}M_t^s(F_i)(x,y)> n\bigg\}$$ are
measurable and decreasing with $n$. If the conclusion of Theorem 5
was false then we could find a measurable set $A\subset (-T/2, T/2)$
with positive measure and a positive number $\delta$ such that for
each $s\in A$ and for each $n\in \N$ we would have
\begin{equation}
H_n(s) > \delta
\end{equation}

We can observe that the set $A$ can be written as
$$A = \bigcap_{n=1}^{\infty}\bigcup_{i=1}^{\infty}\bigg\{s\in(-T/2, T/2):
\mu\bigg\{(x,y)\in \mathcal{D}_N: M_*^s(F_i)(x,y)>
n\bigg\}>\delta\bigg\}.$$

For each $n$ the set $$\bigcup_{i=1}^{\infty}\bigg\{s\in(-T/2, T/2):
\mu\bigg\{(x,y)\in \mathcal{D}_N: M_*^s(F_i)(x,y)>
n\bigg\}>\delta\bigg\},$$ being open by Lemma 5 , it is a countable
union of disjoint open intervals. Therefore the collection (with
$n$) of all these intervals is countable. Because of the decreasing
nature of the sets $$\bigcup_{i=1}^{\infty}\bigg\{s\in(-T/2, T/2):
\mu\bigg\{(x,y)\in \mathcal{D}_N: M_*^s(F_i)(x,y)>
n\bigg\}>\delta\bigg\}$$ with $n,$ the intervals obtained at stage
$k+1$ are included in those corresponding to stage $k.$


\vskip1ex

Our goal is to find a more appropriate countable covering of $A.$
First we can pick an integer $N_1$ large enough and an increasing
sequence of integers $(N_k)_{k>1}$ such that the following
conditions are satisfied.
\begin{enumerate}
\item $$A\subset V_{N_1}=\bigcup_{i=1}^{\infty}\bigg\{s\in (-T/2, T/2):
\mu\bigg\{(x,y)\in \mathcal{D}_N: M_*^s(F_i)(x,y)>
N_1\bigg\}>\delta\bigg\}.$$
\item $$m \bigg\{\bigcup_{i=1}^{\infty}\bigg\{s\in(-T/2, T/2):
\mu\bigg\{(x,y)\in \mathcal{D}_N: M_*^s(F_i)(x,y)>
N_1\bigg\}>\delta\bigg\}\bigg\}\leq 2m(A)$$
\item $\dis \sum_{k=1}^{\infty} \frac{k}{N_k}\leq
(\frac{\delta}{3})^2 m(A)\gamma,$
 where the constant $\gamma$ will
be specified later in order to establish a contradiction.
\end{enumerate}

To start the selection process we pick any $s_1\in A.$ Then there
exists an open interval $I_{1, N_1}\subset V_{N_1}$ that contains
$s_1.$ Then we pick $s_2\in A\cap I_{1, N_1}^c$ and select $I_{1,
N_2}$ containing $s_2.$ By induction we can obtain a countable
collection of open intervals $\dis
\mathcal{J}_1=\bigcup_{k=1}^{\infty} I_{1, N_k}\subset V_{N_1}.$ If
this collection does not cover $A$ then we continue the selection
process by picking $s'\in A\cap \mathcal{J}_1^c$ and an open
interval

\vskip1ex

 $\dis I_{2, N2}\subset V_{N_2}=
\bigcup_{i=1}^{\infty}\bigg\{s\in (-T/2, T/2): \mu\bigg\{(x,y)\in
\mathcal{D}_N: M_*^s(F_i)(x,y)> N_2\bigg\}>\delta\bigg\}\subset
V_{N_1}$ that contains $s'.$ The difference between the collections
$\mathcal{J}_1$ and $\mathcal{J}_2$ is that the first is built with
the sequence $(N_k)_{k\geq 1}$ while the second starting with $N_2$
is built with the sequence $(N_{k+1})_{k\geq 1}$
 Because, as we noticed above, we started with at most countably many open intervals and that at each
 step we picked a different open interval, the selection
 process has to stop after countably many iterations.
 So we obtain after induction at most a countable number of collections $\mathcal{J}_r,$
 $r\in \N,$ that will cover $A$ and will all be contained in $V_{N_1}.$

 \vskip1ex

 We denote the union of these collections of sets by $\dis \mathcal{R} = \bigcup_{r=1}^{\infty}\mathcal{J}_r,.$
 We can observe that with this selection process we have at most one
 interval associated with $N_1$, two with $N_2$ and generally at most $k$ with
 $N_k.$ Now we can use Lemma 6 to extract of this collection of open intervals, disjoint open intervals $G_1$, $G_2,$
 ...$,G_R$ such that

 \vskip1ex

 $(**) \dis  \sum_{h=1}^R m(G_h)> \frac{1}{3} m(A).$

  As all these intervals are disjoint subsets of $V_{N_1}$ we also have
  \vskip1ex

 $\dis(***)\,  \sum_{h=1}^R m(G_h)\leq 2m(A).$

 \vskip1ex
Now we can reach a contradiction. We combine what we obtained so far
to make our choice of $\gamma.$ We have

  \[\begin{aligned}
  &\frac{\delta}{3} m(A)\leq \int_{-T/2}^{T/2} \sum_{h=1}^R
  \mathbf{1}_{G_h}(s)\mu\bigg\{ X\in\mathcal{D}_N;
  M_*^{s}(F_{m_h})(X)> \Gamma_h\bigg\}ds \\
  &\text{for some integers $m_h$ and $\Gamma_h,$} \\
  & \leq \sum_{h=1}^R (m(G_h))^{1/2}\bigg(\int_{-T/2}^{T/2}\bigg(\mu\bigg\{X\in\mathcal{D}_N;
  M_*^{s}(F_{m_h})(X)> \Gamma_h\bigg\}\bigg)^2 ds\bigg)^{1/2} \\
  &\text{by Cauchy Schwartz's inequality,} \\
  & \leq \sum_{h=1}^R (m(G_h))^{1/2}(\mu(\mathcal{D}_N))^{1/2}
  \bigg(\int_{-T/2}^{T/2}\bigg(\mu\bigg\{X\in\mathcal{D}_N;
  M_*^{s}(F_{m_h})(X)> \Gamma_h\bigg\}\bigg)ds\bigg)^{1/2} \\
  & \leq \big(\sum_{h=1}^R m(G_h)\big)^{1/2}(\mu(\mathcal{D}_N))^{1/2}\bigg(\sum_{h=1}^R
  \int_{-T/2}^{T/2}\bigg(\mu\bigg\{X\in\mathcal{D}_N;
  M_*^{s}(F_{m_h})(X)> \Gamma_h\bigg\}\bigg)ds\bigg)^{1/2} \\
  &\text{by Cauchy Schwartz's inequality,}\\
  & \leq (2m(A))^{1/2}\mu(\mathcal{D}_N)^{1/2}T^{1/2}\bigg(\sum_{h=1}^R \frac{4\pi^2 (1+ TK)^2}{(1-TK)^2}\frac{1}{M_h}
   \bigg)^{1/2}\\
  &\text{by using (***), Theorem 3 and the fact that $\|F_{m_h}\|_1\leq 1$} \\
  \end{aligned}
  \]

Therefore we have
\[\begin{aligned}
&\frac{\delta}{3} (m(A))^{1/2}\leq
T^{1/2}\mu(\mathcal{D}_N)^{1/2}(\frac{8\pi^2 (1+
TK)^2}{(1-TK)^2})^{1/2}\bigg(\sum_{h=1}^R \frac{1}{\Gamma_h}\bigg)^{1/2} \\
&\leq T^{1/2}\mu(\mathcal{D}_N)^{1/2}(\frac{8\pi^2 (1+
TK)^2}{(1-TK)^2})^{1/2}\bigg(\sum_{h=1}^{\infty}
\frac{h}{N_h}\bigg)^{1/2} \\
& \text{ because we had for each $k$ at most $k$ intervals corresponding to $N_k$}\\
&< T^{1/2}\mu(\mathcal{D}_N)^{1/2}(\frac{8\pi^2 (1+
TK)^2}{(1-TK)^2})^{1/2}\frac{\delta}{3}(m(A))^{1/2}\gamma^{1/2}\\
&\text{by using (3).}
\end{aligned}
\]
To establish a contradiction it is enough now to pick
$$\gamma< \frac{1}{T\mu(\mathcal{D}_N)\frac{8\pi^2 (1+
TK)^2}{(1-TK)^2}}$$
 choice that we could have made independently of the selection
 process.
 This ends the proof of Theorem 5.

 \vskip1ex
 Because of Lemma 1 Theorem 5 can be reformulated in the following
 way

\begin{corollary}
Let $v$ be a unit Lipschitz vector field with Lipschitz constant
$K= 1/T$. Then there exists a set $\mathcal{T}\subset [-T/2, T/2]$
of measure $T$ such that for each $s\in \mathcal{T}$, for all
$F\in L^p(\R^2)$, $1\leq p<\infty$, the averages
$$\frac{1}{2t}\int_{-t}^t F[(x,y) + (\beta +s)v(x,y))d\beta$$
converge a.e. to $F[(x,y)+ sv(x,y)].$

\end{corollary}

\noindent{\bf Proof of Theorem 4}

 \vskip1ex

 Theorem 5 provides us with a set $\mathcal{T}$ of measure $T$
 such that for each $s\in \mathcal{T}$ we have
 $$\lim_{n}\sup_{\|F\|_1\leq 1}\mu\bigg\{X\in \mathcal{D}_N:
 M_*^s(F)(X)> n \bigg\} = 0.$$

We can conclude, by using similar arguments as those displayed in
\cite{Folland}, that for each $s\in \mathcal{T}$ the set of
functions in $L^1(\R^2)$ for which the pointwise convergence holds
on $\mathcal{D}_N$ is closed in $L^1$. As we obviously have the
dense set of continuous functions for which the differentiation
holds we have proved Theorem 4 for values of $X\in \mathcal{D}_N.$
The general case follows by letting $N$ tend to infinity.

\subsection{ A maximal inequality for the unit vector field $v\circ S_s^{-1}$.}

 First we want to refine the proof of Theorem 5 in order to evaluate the rate of
 convergence to zero (with $n$) of the maximal function
 $\dis \sup_{\|F\|_1\leq 1}\mu\bigg\{(x,y)\in
\mathcal{D}_N: \sup_{0<t<T/2}M_t^s(F)(x,y)> n\bigg\}.$
\begin{theorem}
For each $0<\alpha< 1/2,$ for a.e. $s$ in a set of measure $T$ in
$[-T/2, T/2]$ we have
$$\lim_n n^{\alpha}\sup_{\|F\|_1\leq 1}\mu\bigg\{(x,y)\in \mathcal{D}_N:
\sup_{0<t<T/2}M_t^s(F)(x,y)> n\bigg\}=0.$$
\end{theorem}
\begin{proof}
As in Theorem 5 we argue by contradiction. Instead of the functions
$H_n$ we use this time the functions
$$O_n: s\in [-T/2, T/2]\rightarrow O_n(s) = n^{\alpha}\sup_{\|F\|_1\leq 1}\mu\bigg\{(x,y)\in
\mathcal{D}_N: \sup_{0<t<T/2}M_t^s(F)(x,y)> n\bigg\}$$
 which for each $s$ are equal to
 $$n^{\alpha}\sup_{F_i\in \mathcal{E}}\mu\bigg\{(x,y)\in
\mathcal{D}_N: \sup_{0<t<T/2}M_t^s(F)(x,y)> n\bigg\}$$ (by Lemma 4.)
 The set that replaces $A$ is
 $$B= \bigg\{s \in (-T/2, T/2): \limsup_n O_n(s)>\delta\bigg\}$$
By Lemma 5, for each positive integer $L$ the set
$$ W_L= \bigcap_{n=1}^L\bigcup_{j\geq n}\bigcup_{i=1}^{\infty}\bigg\{s\in (-T/2, T/2):
j^{\alpha}\mu\bigg\{(x,y)\in \mathcal{D}_N:
M_*^s(F_i)(x,y)>j\bigg\}>\delta \bigg\}$$  as a finite intersection
of open sets is a countable union of disjoint open intervals. By
taking the collection (with $L$) of all these open intervals we
obtain a countable number of such intervals. As before the intervals
obtained at stage $L+1$ are subsets of those corresponding to stage
$L.$ Having a countable number of intervals we proceed with an
increasing sequence of integers $N_k$ such that
\begin{equation}
\sum_{k=1}^{\infty}\frac{k}{N_k^{1-2\alpha}}\leq
(\frac{\delta}{3})^2m(B)\gamma'
\end{equation}
We can start the selection process with the additional conditions
\begin{equation}
 B\subset W_{N_1}= \bigcap_{n=1}^{N_1}\bigcup_{j\geq n}\bigcup_{i=1}^{\infty}
\bigg\{s\in (-T/2, T/2): \mu\bigg\{(x,y)\in \mathcal{D}_N:
M_*^s(F_i)(x,y)> j\bigg\}>\delta\bigg\}
\end{equation}
\begin{equation}
 m (W_{N_1})\leq 2m(B)
\end{equation}
As before we select by induction a covering $\dis \mathcal{R}'=
\bigcup_{r'=1}^{\infty}\mathcal{I}_r'$ of $B$ by open intervals,
subsets of those composing $\dis W_{N_1}.$ Furthermore in this
entire collection $\dis \mathcal{R}'$ of open intervals we have at
most one associated with $N_1$, two with $N_2$ and more generally at
most $k$ with $N_k.$ Next we use Lemma 6 to extract of this
collection , disjoint open intervals, $G_1',G_2',...,G_{R}'$ such
that \vskip1ex $\dis (+++)\,  \frac{1}{3} m(B)<\sum_{h=1}^R m(G_h')
\leq 2m(B),$ as these intervals are disjoint subsets of $W_{N_1}.$
To establish the contradiction we will choose later $\gamma'$
appropriately.
 We have
\[
\begin{aligned}
&\frac{\delta}{3} m(B)\leq \int_{-T/2}^{T/2} \sum_{h=1}^R
  \mathbf{1}_{G_h'}(s)(\Gamma_h')^{\alpha}\mu\bigg\{ X\in\mathcal{D}_N;
  M_*^{s}(F_{m_h'})(X)> \Gamma_h'\bigg\}ds \\
  &\text{for some integers $m_h'$ and $\Gamma_h'.$} \\
\end{aligned}
\]
We can use Cauchy Schwarz's inequality to dominate this last term.
\[\begin{aligned}
  & \leq \sum_{h=1}^R (m(G_h'))^{1/2}(\Gamma_h')^{\alpha}\bigg(\int_{-T/2}^{T/2}\bigg(\mu\bigg\{X\in\mathcal{D}_N;
  M_*^{s}(F_{m_h})(X)> \Gamma_h'\bigg\}\bigg)^2 ds\bigg)^{1/2} \\
  & \leq \sum_{h=1}^R (m(G_h))^{1/2}(\Gamma_h')^{\alpha}(\mu(\mathcal{D}_N))^{1/2}
  \bigg(\int_{-T/2}^{T/2}\bigg(\mu\bigg\{X\in\mathcal{D}_N;
  M_*^{s}(F_{m_h'})(X)> \Gamma_h'\bigg\}\bigg)ds\bigg)^{1/2} \\
  & \leq \big(\sum_{h=1}^R m(G_h')\big)^{1/2}(\mu(\mathcal{D}_N))^{1/2}\bigg(\sum_{h=1}^R
  (\Gamma_h')^{2\alpha}\int_{-T/2}^{T/2}\bigg(\mu\bigg\{X\in\mathcal{D}_N;
  M_*^{s}(F_{m_h'})(X)> \Gamma_h'\bigg\}\bigg)ds\bigg)^{1/2} \\
  &\text{by Cauchy Schwartz's inequality,}\\
  & \leq (2m(B))^{1/2}\mu(\mathcal{D}_N)^{1/2}T^{1/2}\bigg(\sum_{h=1}^R \frac{4\pi^2 (1+ TK)^2}{(1-TK)^2}
  \frac{(\Gamma_h')^{2\alpha}}{\Gamma_h'}
   \bigg)^{1/2}\\
   &\text{by using (+++) and Theorem 3} \\
 &=(2m(B))^{1/2}\mu(\mathcal{D}_N)^{1/2}T^{1/2}\bigg(\sum_{h=1}^R \frac{4\pi^2 (1+ TK)^2}{(1-TK)^2}
  \frac{1}{(\Gamma_h')^{1-2\alpha}}\bigg)^{1/2}. \\
&\text{Hence we have}\\
&\frac{\delta}{3}( m(B))^{1/2}\leq
T^{1/2}\mu(\mathcal{D}_N)^{1/2}(\frac{8\pi^2 (1+
TK)^2}{(1-TK)^2})^{1/2}\bigg(\sum_{h=1}^R \frac{1}{\Gamma_h'^{1-2\alpha}}\bigg)^{1/2} \\
&\leq T^{1/2}\mu(\mathcal{D}_N)^{1/2}(\frac{8\pi^2 (1+
TK)^2}{(1-TK)^2})^{1/2}\bigg(\sum_{h=1}^{\infty}
\frac{h}{N_h^{1-2\alpha}}\bigg)^{1/2} \\
& \text{ because we had for each $k$ at most $k$ intervals corresponding to $N_k$}\\
&< T^{1/2}\mu(\mathcal{D}_N)^{1/2}(\frac{8\pi^2 (1+
TK)^2}{(1-TK)^2})^{1/2}\frac{\delta}{3}(m(B))^{1/2}\gamma'^{1/2}\\
&\text{by using (4).}
\end{aligned}
\]

The contradiction is obtained for $$
\gamma'<\frac{1}{T\mu(\mathcal{D}_N)\frac{8\pi^2 (1+
TK)^2}{(1-TK)^2}}.$$

\end{proof}
 The following result can be derived from Theorem 5.
\vskip1ex
\begin{theorem}
For each $0<\alpha<1/2$ there exists a function $C_{\alpha}$ almost
everywhere finite on $[-T/2, T/2]$ such that for all $\lambda >1,$
for all $F\in L^1$ such that $\lambda>\|F\|_1$ we have
\begin{equation}
\mu\bigg\{X\in \mathcal{D}_N: M_*^s(F)(X)>\lambda\bigg\}\leq
2^{\alpha}C_{\alpha}(s)\bigg(\frac{\|F\|_1}{\lambda}\bigg)^{\alpha}
\end{equation}
\end{theorem}
\begin{proof}
 For a fixed $\alpha$ we denote by $C_{n,\alpha}(s)$ the a.e. finite
 function $$n^{\alpha}\sup_{\|F\|_1\leq 1}\mu\bigg\{(x,y)\in \mathcal{D}_N:
\sup_{0<t<T/2}M_t^s(F)(x,y)> n\bigg\}.$$ By Theorem 6 we have $\dis
\lim_n C_{n,\alpha}(s) = 0 $ for a.e. $s\in [-T/2, T/2].$ Hence the
function $\dis C_{\alpha}: C_{\alpha}(s)= \sup_{n} C_{n,\alpha}(s)$
is a.e. finite on $[-T/2, T/2].$ Furthermore for each $\lambda >1$
we have
\[\begin{aligned}
&\sup_{\|F\|_1\leq 1}\mu\bigg\{(x,y)\in \mathcal{D}_N:
\sup_{0<t<T/2}M_t^s(F)(x,y)> \lambda\bigg\}\\
&\leq \sup_{\|F\|_1\leq 1}\mu\bigg\{(x,y)\in \mathcal{D}_N:
\sup_{0<t<T/2}M_t^s(F)(x,y)> [\lambda]\bigg\}\\
&\leq
\frac{C_{\alpha}(s)}{[\lambda]^{\alpha}}\\
\end{aligned}
\]
 Therefore by using the inequality $\lambda <2 [\lambda]$ we have
  $$ \sup_{\|F\|_1\leq 1}\mu\bigg\{X\in \mathcal{D}_N:
M_*^s(F)(X)>\lambda\bigg\} \leq
\frac{2^{\alpha}C_{\alpha}(s)}{\lambda^{\alpha}}$$
 Now by taking the functions $F/\|F\|_1$ and making the change
$\lambda \|F\|_1 = t$ we can derive (7).

\end{proof}
\section{Differentiation in $R^n$}
The results obtained in the previous section can be extended without
difficulty to $\R^n.$ In fact the only part where we use the fact
that we were in $\R^2$ is when we proved Lemma 1. This appears in
the constant of this Lemma. We only state the lemma that would
replace Lemma 1 in $\R^n.$ We consider a Lipschitz unit vector field
$v$ on $R^n$ with constant $K$ and simply denote by $\dis M_t(F)(X)$
the averages
$$\frac{1}{2t}\int_{-t}^t F[X+ \beta v(X)]d\beta$$ where $X = (x_1,x_2, ...,x_n).$
We use the same notation for $S_s(X) = X + \beta v(X)$ a map from
$R^n$ to $R^n.$  We denote by $\mu_n$ Lebesgue measure on $R^n.$
\begin{lemma}
For all $|s| \leq T$ where $T<1/K$ the maps $S_s$ are one to one and
onto. Furthermore there exist constants $c_n$ and $C_n$ depending
only on $n$, $K$ and $T$ such that for all measurable sets $A$ in
$R^n$ we have
$$c_n\mu_n(S_s(A))\leq \mu_n(A) \leq C_n\mu_n(S_s(A))$$
\end{lemma}
\begin{proof}
 The invertible character of the maps $S_s$ for small $s$ can be
 established in the same way.  The inequalities
$$c_n\mu_n(S_s(A))\leq \mu_n(A) \leq C_n\mu_n(S_s(A))$$
follow from the inequalities $\|Z_1- Z_2\|\leq (1+ |s|K)\|X_1-
X_2\|,$ and
 $\|X_1-X_2\|\leq \frac{1}{1-|s|K} \|S_s(X_1)-S_s(X_2)\|$
 where $Z_1 = S_s(X_1)$ and $Z_2= S_s(X_2).$
\end{proof}
 The maximal inequality that replaces Theorem 3 is the following
 \begin{theorem}
 Let $K$ be the Lipschitz constant for the unit vector field $v.$ Then for each $T$, $0<T< 1/K,$
 there exists a constant $\mathcal{C}_n$ such that for all $\lambda>0$
\[
\begin{aligned}
&\frac{1}{T}\int_{-T/2}^{T/2}\mu_n\bigg\{X\in R^n: \sup_{0<t\leq
T/2}\frac{1}{2t}\int_{-t}^t |F[ X + \beta
v(S_s^{-1}(X))]|d\beta >\lambda\bigg\}dm(s)\\
&\leq \frac{\mathcal{C}_n}{\lambda}\int_{\R^n} |F(X)|d\mu_n.
\end{aligned}
\]
 where $m$ denotes Lebesgue measure on $[-T/2, T/2].$
\end{theorem}
The proof is identical to the one given for Theorem 3 so we skip it.
From this maximal inequality the reasoning is identical. The only
difference is the constant $\mathcal{C}_n$ that depends on $T,$ $K$
and $n.$  From that point on by using the same path one can extend
Theorem 4 and Theorem 5 to the case of $\R^n.$


\begin{thebibliography}{99}

\bibitem{Bourgain}{\bf {J. Bourgain}}:``A remark on the maximal
function associated to an analytic vector field'', \textit{Proc. Of
the Special Year in Modern Analysis,} London Math. Soc. Lect. note
Series, vol 137, 111-132, 1989.
\bibitem{DeGuzman}{\bf M. De Guzman}:``Real variable methods in Fourier Analysis'',
\textit{Math. Study,} 46, North Holland, 1981.
\bibitem{Folland}{\bf G. Folland}:``Real Analysis",
\textit{J. Wiley and Sons,} 1984.
\bibitem{Katz}{\bf {N.H. Katz}}:``A partial result on Lipschitz
differentiation'', \textit{ Harmonic Analysis at Mount Holyoke,}
Contemporary Math. , vol 320, 217-224, 2003.
\bibitem{LaceyLi}{\bf {M. Lacey and X. Li}}:``On the Hilbert Transform and $C^{1+\epsilon}$ Families of
Lines'', \textit{Annals of Math.}, to appear.
\bibitem{Stein}{\bf E.M. Stein}:``Real variable methods,
orthogonality, and oscillatory integrals'', \textit{Princeton Math.
Series,} 43, 1993.
\bibitem{Wiener}{\bf N. Wiener } :``The ergodic theorem'',
\textit{ Duke Math. J.}, 5, 1-18, 1939.
 \end{thebibliography}
\end{document}